\documentclass{tsp-short}
\newtheorem{thm}{Theorem}[section]
\newtheorem{lemma}{Lemma}[section]

\newtheorem{assertion}{Proposition}[section]
\theoremstyle{remark}
\newtheorem{rem}{Remark}[section]
\theoremstyle{definition}
\newtheorem{define}{Definition}[section]
\newtheorem{ex}{Example}[section]

\usepackage[dvips]{graphicx}

\newcommand{\mmp}{\mathbb{P}}

\newcommand{\od}{\overset{d}{=}}
\newcommand{\bydef}{\overset{def}{=}}

\newcommand{\ind}{\mathbb{\text{1}}}

\newcommand{\dod}{\overset{d}{\to}}
\newcommand{\tp}{\overset{P}{\to}}

\newcommand{\me}{\mathbb{E}}
\newcommand{\mr}{\mathbb{R}}
\newcommand{\mn}{\mathbb{N}}
\newcommand{\mno}{\mathbb{N}_0}

\newcommand{\lin}{\underset{n\to\infty}{\lim}}

\newcommand{\lix}{\underset{x\to\infty}{\lim}}

\renewcommand{\labelenumi}{\theenumi)}

\begin{document}
\title{On the asymptotics of moments of linear random recurrences}

\author{Alexander Marynych}
\address{Address of Alexander Marynych}
\curraddr{Faculty of Cybernetics, National T. Shevchenko University of Kiev, 01033 Kiev, Ukraine}
\email{marynych@unicyb.kiev.ua}

\date{12/05/2010}

\begin{abstract}
\noindent
We propose a new method of analyzing the asymptotics of moments of
certain linear random recurrences which is based on the technique
of iterative functions. By using the method, we show that the
moments of the number of collisions and the absorption time in the
Poisson-Dirichlet coalescent behave like the powers of the "log
star" function which grows slower than any iteration of the
logarithm, and thereby prove a weak law of large numbers. Finally,
we discuss merits and limitations of the method and give several
examples related to beta coalescents, recursive algorithms and
random trees.
\end{abstract}

\maketitle
\section{Introduction and main result.}
\paragraph{}
A {\it linear random recurrence} is a sequence of random variables
$\{X_n,n\in\mn\}$ which satisfies the distributional equality
\begin{eqnarray}\label{rec_main}
X_1=0, \ \ \ \ X_n\od V_n+\sum_{r=1}^{K}A_r(n)X^{(r)}_{I^n_{(r)}},
\ \ n\geq 2,
\end{eqnarray}
where $X_n$ is some parameter of a problem of size $n$, which
splits into $K\geq 1$ subproblems of random sizes
$I^n_{(r)}\in\{1,\ldots,n\}$. For every $r=1,\ldots,K$, the
sequence $\{X_k^{(r)},k\in\mn\}$ which corresponds to the
contribution of subgroup $r$ is a distributional copy of
$\{X_k,k\in\mn\}$, $V_n$ is a random toll term, and $A_r(n)>0$ is
a random weight of subgroup $r$. It is assumed that
$\{(I^n_{(1)},\ldots,I^n_{(K)},A_1(n),\ldots,A_K(n),V_n),n\geq
2\}$, $\{X_n^{(1)},n\in\mn\},\ldots,\{X_n^{(K)},n\in\mn\}$ are
independent.

Random  recurrences \eqref{rec_main}, often in a simplified form
with $K=1$, arise in diverse areas of applied probability such as
random regenerative structures \cite{gnedin_sieve, GINR}, random
trees \cite{Drmota,DIMR2,Pan1,Pan2}, coalescent theory
\cite{DIMR1,GneIksMoe,GneYak, IMM}, absorption times in
non-increasing Markov chains \cite{Haas, Cutsem}, recursive
algorithms \cite{KnutGrin,NeiRus,RosQu,Roesler99}, random walks
with barrier \cite{IksMoe2,IksNegad}, to name but a few.

The first step of asymptotic analysis of recurrences
\eqref{rec_main} is to find the asymptotics of moments
$\mathbb{E}X_n^k$ and central moments
$\mathbb{E}(X_n-\mathbb{E}X_n)^k$, as $n\to\infty$. This problem
reduces to studying the recurrent equations of the form
\begin{eqnarray}\label{rec_mean}
a_1=0,\;\;a_n=b_n+\sum_{k=1}^{n-1}c_{nk}a_k,\;n\geq 2,
\end{eqnarray}
where $\{b_n,n\in\mn\}$ and $\{c_{nk},n\in\mn,k<n\}$ are given
numeric sequences. The purpose of the present paper is to propose
a new method of obtaining the first-order asymptotics of solutions
to (\ref{rec_mean}), as $n\to\infty$.

Although the asymptotic analysis of recurrences (\ref{rec_mean})
is a hard analytic problem, some more or less effective methods
have been elaborated to date. Evidently the most popular existing
approach is the {\it method of singular analysis of generating
functions} \cite{Drmota,FlajSedg}. The method gives a very precise
information on the asymptotic behavior whenever there is a
tractable functional relation between the generating functions of
the sequences involved. The idea of a {\it repertoire method}
proposed in \cite{KnutGrin} can be briefly described as follows.
First we build up a repertoire $\{b^{(\alpha)}_n,\;\alpha\in A\}$,
where $A$ is a finite set, of inhomogeneous terms of
(\ref{rec_mean}) by choosing sequences
$\{a^{(\alpha)}_n,\;n\in\mn\}$ such that the sum in
(\ref{rec_mean}) is tractable. Then we construct the solution
$a_n$ to (\ref{rec_mean}) with inhomogeneous term $b_n$ as a
linear combination of solutions $a^{(\alpha)}_n,\;\alpha\in A$.
Finally, we mention a method proposed in \cite{Bruhn} and further
developed in \cite{Roesler99} which is based on the harmonic
analysis and potential theory.

%Motivated by the question asked by Martin M\"{o}hle about the
%asymptotics of moments of the number of collisions in the
%Poisson-Dirichlet coalescent we tried unsuccessfully to apply
%known to us existing methid

%we choose $a_n$ such that the sum in (\ref{rec_mean}) is
%tractable, then we calculate corresponding inhomogeneous term
%$b_n$. Having big enough repertoire of defferent additive terms we
%can construct the solution of (\ref{rec_mean}) as a linear
%combination of these terms.

Initially, our method which is based on the technique of iterative
functions was developed in order to find the asymptotics of
moments of the number of collisions $X_n$ and the absorption time
$T_n$ in the Poisson-Dirichlet coalescent (this problem was raised
by Martin M\"{o}hle in the early 2008). Below we prove that both
$\mathbb{E}X_n^k$ and $\mathbb{E}T_n^k$, $k\in\mn$, behave like
the powers of the "log star" function which grows slower than any
iteration of the logarithm\footnote{The result for $\mathbb{E}X_n$
was conjectured by M.~M\"{o}hle.}. This somewhat exotic behavior
of the moments partially explains the fact that we have not been
able to apply either of previously known (to us) methods to tackle
the problem. To our knowledge, the "log star" asymptotics arises
not often. In particular, we are only aware of two applied models
which exhibit such a behavior: (a) the number of distinguishable
alleles according to the Ohta-Kimura model of neutral mutation
\cite{Kesten}, and (b) the average complexity of Delaunay
triangulation of the Euclidian minimum spanning tree
\cite{Devillers}. The number of collisions and the absorption time
in the Poisson-Dirichlet coalescent are interesting, yet
particular patterns of recurrence \eqref{rec_main}. Thus after
having settled the original problem concerning the
Poisson-Dirichlet coalescent, our method was subsequently extended
to cover many other linear recurrences \eqref{rec_main}.

In this paper, unless stated the contrary, we tacitly suppose that
$b_n\geq 0$ and, hence, $a_n\geq 0$. However, a perusal of the
proofs given below reveals that we could have assumed that $b_n$
is only non-negative or non-positive for large enough $n$. Under
this last assumption, formulations of results would get cumbersome
which has forced us to keep less generality but more transparency.

Our method can be summarized in the following
\begin{center}
\begin{LARGE}
Algorithm
\end{LARGE}
\end{center}
\begin{enumerate}
\item Using, for example, the method described in \cite{Roesler99}, obtain the recurrence with weights reduced to probabilities.
As a result, we obtain the recurrence of the form
$$
A_1=0,\;\;A_n=B_n+\sum_{k=1}^{n-1}p_{nk}A_k,
$$
where $\sum_{k=1}^{n-1}p_{nk}=1$ for all $n\geq 2$ and $B_n\geq
0$. Let $I_n$ be a random variable with distribution
$\mmp\{I_n=k\}=p_{nk},n\geq 2,k<n$.
\item Prove the divergence of $A_n$ using, for example, Proposition \ref{is_unbounded} or other methods.
\item Find a continuous, strictly increasing and unbounded function $g(x)$ defined on $\mr^{+}$,
such that $g(n)=\mathbb{E}I_n+o(\mathbb{E}I_n)$. Pick an $x_0$ as
defined in (\ref{condition_iter}). Find a continuous function $h(x)$ defined on $\mr^+$ such that $h(n)=B_n$.
\item Find an iterative function $g^{*}$ generated by the quadruple $(h,g,x_0,k)$,
where $k$ is any continuous on $[0,x_0]$ function (see Definition
\ref{iterative_function_definition}).

\item Using, for example, Theorem \ref{f_equiv} find, if possible, an elementary function $f_1$ such that
$\lim_{x\to\infty}\frac{f_1(x)}{g^{*}(x)}=1$, and set $f:=f_1$.
Otherwise, select $k$ such that $g^{*}$ is twice differentiable,
and set $f:=g^{*}$ (see Theorem \ref{main_aproximation}).

\item If $f(\mathbb{E}I_n)-f(g(n))=o(h(n))$ then go to the next step, otherwise go to step 3) and choose asymptotically smaller term $o(\mathbb{E}I_n)$.
\item If $\mathbb{E}f(I_n)-f(\mathbb{E}I_n)=o(h(n))$ then $A_n\sim f(n)$.
      If $\mathbb{E}f(I_n)-f(\mathbb{E}I_n)\sim ch(n)$ then $A_n\sim (1-c)^{-1}f(n)$ (see Theorem \ref{main_theorem_1}).
\end{enumerate}

We mention, in passing, that iterative functions have already been
used in the context of divide-and-conquer paradigm \cite{Karp}.
The cited paper is concerned with stochastic processes
$\{T(x),x\in\mr^{+}\}$ whose marginal distributions are given by
the equality
$$
T(x)\od a(x)+T'(t(x)),\;\;x\in\mr^{+},
$$
where $a(\cdot)$ is a non-negative (deterministic) function and
$t(\cdot)$ is a random variable taking values in $[0,\cdot]$ which
is independent of $\{T'(x),x\in\mr^{+}\}$, an independent copy of
$\{T(x),x\in\mr^{+}\}$.

The structure of the paper is as follows. Section 2 introduces
iterative functions and investigates their properties. Section 3
carefully describes the algorithm of our new method. Theorem
\ref{theorem_general} and Theorem \ref{main_theorem_1} which are
the main results of the section prove the validity of the
algorithm. Section 4 is devoted to applications and also discusses
"ins and outs" of the method. The paper closes with the Appendix
which collects proofs of some technical results concerning the
iterative functions and properties of recurrences
\eqref{rec_mean}.

Throughout the paper the notation $r(\cdot)\sim s(\cdot)$ means that
$r(\cdot)/s(\cdot)\to 1$, as the argument tends to infinity, 
$C^{(m)}(B)$ denotes the space of functions which are $m$-times
differentiable on the set $B$. If $B=[a,\infty)$ then the
derivatives at point $a$ are assumed to be the right derivatives.
Also we use notation
$$
r^{\circ (0)}(x)\bydef x,\;\; r^{\circ (k)}(x)\bydef r(r^{\circ
(k-1)}(x)),\;\;k\in\mn.
$$
Finally, we recall the standard notation $\lfloor x\rfloor=
\sup\{k\in\mathbb{Z}:k\leq x\}$ and $\lceil x\rceil=
\inf\{k\in\mathbb{Z}:k\geq x\}$ for the floor and ceiling
function, respectively.

\section{Iterative functions.}\label{2}

In this section iterative functions are defined and some basic
properties of these functions are given. We start with a formal
definition.

\begin{define}\label{iterative_function_definition}
Suppose that the function $g:\mr^+ \to\mr^+$ is increasing,
unbounded and continuous, and satisfies the following condition:
for some $x_0>0$ and every $x_1>x_0$ there exists
$\varepsilon_{x_1}>0$ such that
\begin{eqnarray}\label{condition_iter}
x-g(x)>\varepsilon_{x_1} \ \ \text{for all} \ \  x\in(x_0,x_1).
\end{eqnarray}
Assume that the functions $h:\mr^+ \to\mr^+$ and $k:[0,x_0]\to
\mr$ are continuous and define the function $g^*:\mr^{+}\to\mr$ by
the following equality
\begin{eqnarray}\label{sum_representation}
g^*(x)=\sum_{i=1}^{m_0(x)}h(g^{\circ (i-1)}(x))+k(g^{\circ (m_0(x))}(x)),
\end{eqnarray}
where
$$
m_0(x):=\inf\{k\geq 0:g^{\circ (k)}(x)\leq x_0\}.
$$
We call $g^{*}$ the {\it iterative function generated by the
quadruple} $(h,g,x_0,k)$ and denote it by $g^{*}={\rm
Iter}(h,g,x_0,k)$.
\end{define}
Note that technical condition (\ref{condition_iter}) is sufficient
for $m_0(x)$ to be finite for every $x\in\mr^+$. This follows from
the estimate $m_0(x)\leq
\lfloor\frac{x-x_0}{\varepsilon_x}\rfloor+1$, which is implied by
the inequality
$$
x-k\varepsilon_x>g^{\circ (k)}(x), \ \ x>x_0, \ \ k\in\mn,
$$
which in turn can be obtained by induction.
\begin{rem}
From the definition it follows that $g^*$ satisfies the functional equation
\begin{eqnarray}\label{iter_def}
g^*(x)
 & = & h(x)+g^*(g(x)),\; x>x_0,
\end{eqnarray}
with initial condition
\begin{eqnarray*}\label{init_condition}
g^*(x)
 & = & k(x),\qquad\qquad\;\;\; x\leq x_0.
\end{eqnarray*}
\end{rem}
Below are some examples of iterative functions
\paragraph{}
\textit{Example 1}. Let $h(x)\equiv 1$, $g(x)=\alpha x$,
$\alpha\in [0,1)$, $x_0=1$, $k(x)\equiv 0$. Then $g^*(x)=1+g^*(\alpha x),\;\;x>1$,
or $$g^*(x)=\lceil \log_{\frac{1}{\alpha}} x\rceil,\;\;x>1.$$

\paragraph{}
\textit{Example 2}. Let $h(x)\equiv 1$, $g(x)=\log x$, $x_0=1$,
$k(x)\equiv 0$. Then $g^*(x)=1+g^*(\log x ),\;\;x>1$, or
$$g^*(x)=\log^* x,$$ the {\it log-star function} which is arguably the best known non-trivial
iterative function. It is clear that ${\rm Iter} (1,g,x_0,
0)=m_0(x)$. In particular, this equality holds for the log-star
function.
\paragraph{}

If $h(x_0)\neq 0$ then the iterative functions ${\rm Iter}
(h,g,x_0,0)$ are piecewise continuous. We however prefer to work
with smooth iterative functions which was the main reason for
introducing functions $k$ in Definition
\ref{iterative_function_definition}. It turns out that ${\rm
Iter}(h,g,x_0, 0)$ and ${\rm Iter} (h,g,x_0, k)$ have the same
asymptotics, and an appropriate choice of $k$ makes ${\rm Iter}
(h,g,x_0, k)$ smooth enough. Below we formalize this statement and
also describe how the mentioned smoothness can be obtained by the
choice of $k$.

%It had already been said in the introduction that the log-star
%function arises in some applications.  The class of iterative
%function generated by quadruple $(1,g,x_0,0)$ (which
%covers the log-star function) has good interpretation. In this
%case $g^{*}(x)\equiv m_0(x)$, which means that $g^{*}$ is equal to
%the number of times (steps) one should apply the function $g$ to
%the number $x$  to lower the level $x=x_0$.
%
%Apparently the most natural way of generalizing this class is to
%introduce the counting function $h$, which weights steps. On the
%other hand the primary reason for introducing the initial function
%$k$ in the definition was to provide iterative function with good
%analytical properties, particularly with the smoothness. In this
%paper we are interested mainly in the asymptotic behavior of
%iterative function and it should be mentioned that function $k$
%(in contrast to $h$) does not influence asymptotic properties of
%iterative function. Below we formalize this statement and also
%describe how the mentioned smoothness can be provided by the
%choice of the function $k$.

Introduce the equivalence relation $\approx$ on the set of iterative functions by the rule
$$
g_1^{*}\approx g_2^{*}\Longleftrightarrow g_1^{*}={\rm
Iter}(h,g,x_0,k_1),\;\;g_2^{*}={\rm Iter}(h,g,x_0,k_2).
$$
This relation induces partitioning the set of iterative functions
into the classes of equivalence.
\begin{define}\label{iter_def_by_triple}
The equivalence class
$$\mathcal{F}:=\{F={\rm Iter}(h,g,x_0,k),k\in C[0,x_0]\}$$ is called the \textit{iterative function generated by the triple} $(h,g,x_0)$.
When it does not lead to ambiguity, we call an \textit{iterative
function generated by the triple} $(h,g,x_0)$ an arbitrary element
of this class.
\end{define}
Since $|g_1^{*}(x)-g_2^{*}(x)|$ is bounded on $\mr^{+}$, for any
$g_1^\ast, g_2^\ast\in \mathcal{F}$, all iterative function in the
same equivalence class are asymptotically equivalent (provided
they diverge).
\begin{define}
An \textit{$m$-times differentiable modification} of iterative
function $g^{*}$ is an arbitrary iterative function $\hat{g}^{*}$
such that $\hat{g}^{*}\approx g^{*}$ and $\hat{g}^{*}\in
C^{(m)}[x_0,+\infty)$.
\end{define}

%Recall that the main reason of introducing function $k$ was to
%make iterative function smooth ( in the meaning that it lies in
%the space $C^{(m)}[x_0,+\infty)$). Further we assume that
%functions $g,h$ are differentiable suitable number of times on the
%set  $[x_0,+\infty)$.

Our first result which is a direct consequence of Lemma \ref{lem1}
and Lemma \ref{lem2} given in the Appendix shows that provided $g$
and $h$ are smooth enough one can find a function $k$ such that
the function ${\rm Iter}(h,g,x_0, k)$ is smooth. For a collection
of functions $f_1,\ldots,f_n$ let $W(f_1,\ldots,f_n)$ denote its
wronskian.
\begin{thm}\label{main_aproximation}
Assume that $g,h\in C^{(m)}[x_0,+\infty)$ and that
$$W\Big(x^i-g^i(x),i=0,\ldots,m+1\Big)(x_0)\neq 0.$$ Then there
exists a function $k$ of the form
$$
k(x)=\sum_{i=1}^{m+1}\alpha_i x^i,
$$
such that the iterative function generated by the quadruple
$(h,g,x_0,k)$ is $m$-times differentiable on $[x_0,+\infty)$.
\end{thm}
\begin{rem}
The vector of coefficients
$(\alpha_1,\alpha_2,\ldots,\alpha_{m+1})$ is a solution to the
system of linear equations (see Lemma \ref{lem2}) and can be
calculated explicitly.
\end{rem}

An example of a smoothed iterative function is given next.
\begin{ex}\label{iter_log}
Recall that the log-star function is iterative function generated
by the quadruple $(1,\log x, 1,0)$. A twice differentiable
modification $F$ of the log-star function can be constructed as
follows. According to Lemma \ref{lem2}, the corresponding function
$k$ takes the form
$k(x)=-\frac{2}{13}x^3+\frac{3}{13}x^2+\frac{12}{13}x$. Therefore
$$
F(x)=
\left
\{
\begin{array}{ll}
1+F(\log x),\;\;x>1,\\
-\frac{2}{13}x^3+\frac{3}{13}x^2+\frac{12}{13}x, x\in [0,1].
\end{array}
\right.
$$
Below are depicted the graphs of functions $\log^* x$ and $F(x)$
for $x>0$.
\begin{center}
\includegraphics[angle=-90,scale=0.4]{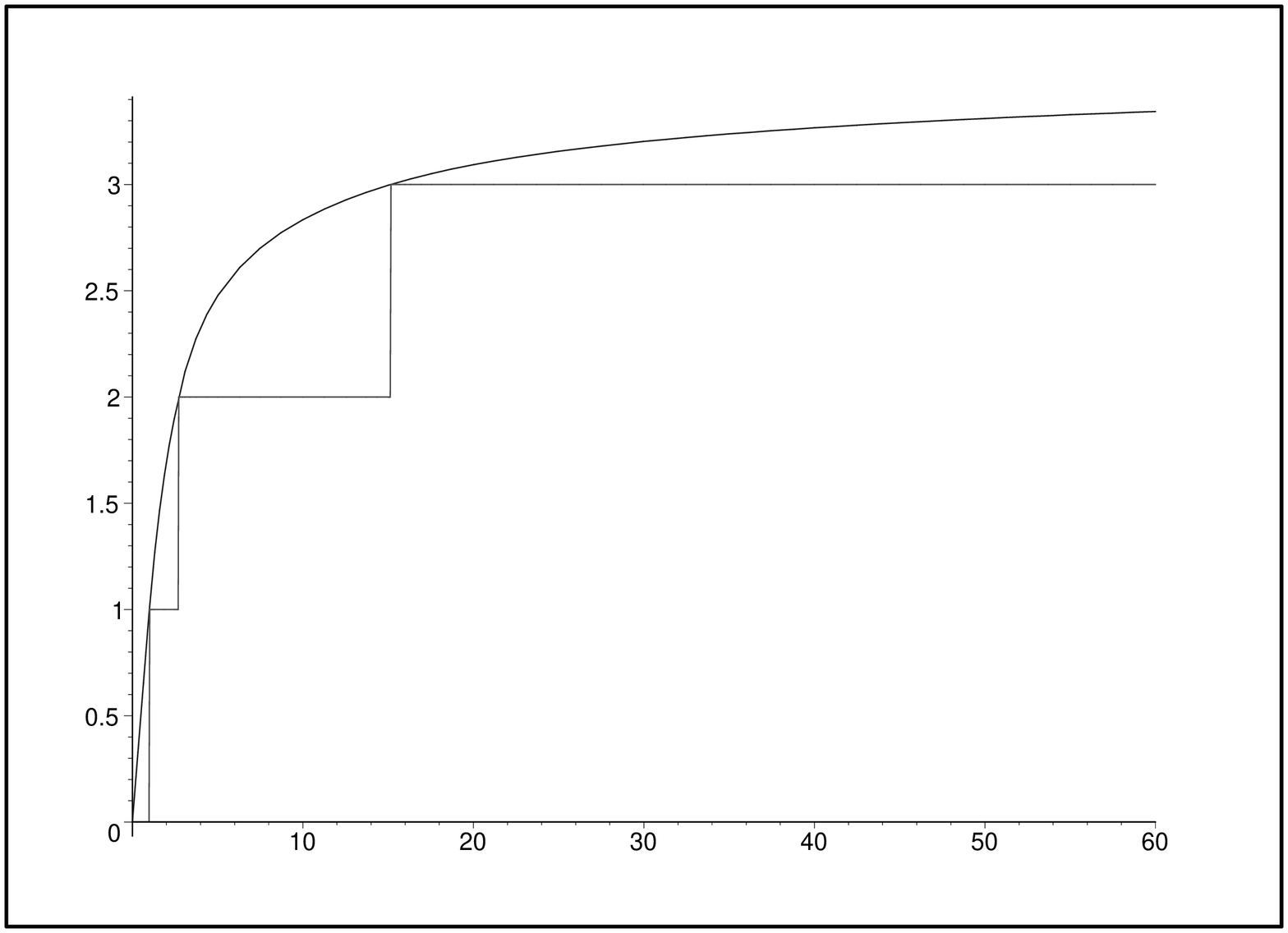}
\end{center}
\end{ex}

\section{Asymptotic behavior of (\ref{rec_mean}).}
While investigating recurrence \eqref{rec_mean}, without loss of
generality, we can assume that for every $n\geq 2$
\begin{equation}\label{3}
\sum_{k=1}^{n-1}c_{nk}=1\text{ and } c_{nk}\geq 0,\;k=1,\ldots,n-1
\end{equation}
(see, for instance, p.~9 in \cite{Roesler99}). In what follows,
recurrences \eqref{rec_mean} with $b_n\geq 0$ which satisfy
\eqref{3} are referred to as {\it recurrences with weights reduced to
probabilities}. If \eqref{3} holds, denote by $I_n$ a random
variable with distribution
$$\mmp\{I_n=k\}=c_{nk}, \ \ k=1,\ldots,n-1.$$ %In what follows, we
%use the functions $g$ and $h$ as defined in Section \ref{2}.
\begin{thm}\label{theorem_general}
Assume that the sequence $\{a_n,n\in\mn\}$ satisfy recurrence
\eqref{rec_mean} with weights reduced to probabilities.
%$$a_1=0,\;\;a_n=b_n+\sum_{k=1}^{n-1}c_{nk}a_k, \ \ n\geq 2.$$
Let $g:\mr^+\to\mr^+$ be a continuous, increasing and unbounded
function such that $$g(n)=\me I_n+o(\me I_n), \ \ n\to\infty,$$
and $h:\mr^+\to\mr^+$ be a continuous function such that
$$h(n)=b_n, \ \ n\geq 2.$$ If
\begin{itemize}
\item $\lin a_n=+\infty$,
\item $g^{*}(\mathbb{E}I_n)-g^{*}(g(n))=o(h(n))$, $n\to\infty$,
\end{itemize}
where $g^{*}$ is an iterative function generated by the triple
$(h,g,x_0)$, then the following implications are true
\begin{eqnarray}\label{conclusion1}
&   &\hspace{-2.7cm}\mathbb{E}g^{*}(I_n)-g^{*}(\mathbb{E}I_n)=o(h(n)), \ \ n\to\infty \Longrightarrow \\
&   & a_n\sim g^{*}(n), \ \ n\to\infty,\nonumber
\end{eqnarray}
\begin{eqnarray}\label{conclusion2}
&   &\hspace{-0.1cm}\mathbb{E}g^{*}(I_n)-g^{*}(\mathbb{E}I_n)\sim dh(n), \ \ n\to\infty, \ \ \text {for some }d<1 \Longrightarrow \\
&   &\hspace{2.6cm}a_n\sim (1-d)^{-1}g^{*}(n), \ \ n\to\infty.
\nonumber
\end{eqnarray}
\begin{proof}
Set $a'_n:=a_n-g^{*}(n)$, $n\in\mn$. The sequence $\{a'_n, n\in\mn\}$
satisfies the recurrence
\begin{equation}\label{1}
a'_1=-g^{*}(1),\;\;a'_n=b_n-g^{*}(n)+\mathbb{E}g^{*}(I_n)+\sum_{k=1}^{n-1}c_{nk}a'_k,
\ \ n\geq 2.
\end{equation}
If $\mathbb{E}g^{*}(I_n)-g^{*}(\mathbb{E}I_n)=o(h(n))$ and
$g^{*}(\mathbb{E}I_n)-g^{*}(g(n))=o(h(n))$ then the inhomogeneous term of
\eqref{1} is $o(h(n))$. Therefore, applying part (II) of Theorem
\ref{seq_equiv} yields $a'_n=o(a_n)$ which implies $a_n\sim g^{*}(n)$.

If $\mathbb{E}g^{*}(I_n)-g^{*}(\mathbb{E}I_n)\sim dh(n)$ for some
$d\in (0,1)$ and $g^{*}(\mathbb{E}I_n)-g^{*}(g(n))=o(h(n))$ then
the inhomogeneous term of \eqref{1} is asymptotically equal to
$dh(n)$. Therefore, applying part (I) of Theorem \ref{seq_equiv}
yields $a'_n\sim da_n$ which implies $a_n\sim (1-d)^{-1}g^{*}(n)$.

Finally, if $\mathbb{E}g^{*}(I_n)-g^{*}(\mathbb{E}I_n)\sim dh(n)$ for some
$d<0$ we can apply part (II) of Theorem \ref{seq_equiv} to the
sequences $\{g^{*}(n)-a_n\}$ and $\{a_n\}$ to conclude that
$g^{*}(n)-a_n\sim -da_n$. The latter is equivalent to
$a_n\sim(1-d)^{-1}g^{*}(n)$. The proof is complete.
\end{proof}
\end{thm}
\begin{thm}\label{main_theorem_1}
Assume that the sequence $\{a_n,n\in\mn\}$ satisfy recurrence
\eqref{rec_mean} with weights reduced to probabilities. Let $g:\mr^+\to\mr^+$
be a twice differentiable, increasing and unbounded function such
that
$$g(n)=\me I_n+o(\me I_n), \ \ n\to\infty,$$ 
and $h:\mr^+\to\mr^+$ be a twice differentiable function such that
$$h(n)=b_n, \ \ n\geq 2.$$
If the following conditions hold
\begin{enumerate}
 \item[(C1)] $\lin a_n=+\infty$
 \item[(C2)] There exists a continuous function $k$ such that the iterative function $F$ generated by the quadruple
 $(h,g,x_0,k)$ is twice differentiable
 \item[(C3)] $F(\mathbb{E}I_n)-F(g(n))=o(h(n))$, $n\to\infty$
 \item[(C4)] There exists $M>0$ such that for all $n\in\mn$ $${\rm Var}\,I_n \leq M \mathbb{E}I_n$$
 \item[(C5)] $\lin \sup_{x\geq {\mathbb E}I_n/2}|F''(x)|\dfrac{{\rm
 Var}\,I_n}{h(n)}=0$
 \item[(C6)] $\lin \dfrac{\sup_{1\leq x\leq n}|F(x)|}{h(n) {\rm Var}\,I_n}= 0$
 \item[(C7)] $\lin \dfrac{F'(\mathbb{E}I_n)}{h(n)}= 0$,
\end{enumerate}
then
\begin{eqnarray*}\label{equiv1}
a_n\sim F(n),\;n\to\infty.
\end{eqnarray*}
\begin{proof}
Since conditions (C1) and (C3) hold, according to implication
\eqref{conclusion1} in Theorem \ref{theorem_general}, it is enough
to show that
\begin{equation*}\label{4}
\alpha_n:=\mathbb{E}F(I_n)-F(\mathbb{E}I_n)=o(h(n)).
\end{equation*}
With $\kappa:=\frac{1}{2M}$ and $A_n:=\{|I_n-\mathbb{E}I_n|>\kappa
{\rm Var}\,I_n\}$ we have $$|\alpha_n|\leq
|\mathbb{E}(F(I_n)-F(\mathbb{E}I_n))\ind_{A_n}|+|\mathbb{E}(F(I_n)-F(\mathbb{E}I_n))\ind_{A_n^c}|=:\beta_n+\gamma_n.$$
An application of Chebyshev's inequality yields
$$
\beta_n \leq 2\sup_{1\leq x\leq n}|F(x)|\mmp(A_n)\leq \frac{2{\rm
Var}\,I_n\sup_{1\leq x\leq n}|F(x)|}{(\kappa {\rm Var}\,I_n)^2},
$$
which is $o(h(n))$ by condition (C6).

Using the Taylor expansion around $\mathbb{E}I_n$ leads to
\begin{eqnarray*}
\gamma_n
& = & \Big|\mathbb{E}\Big(F'(\mathbb{E}I_n)(I_n-\mathbb{E}I_n)+\frac{1}{2}F''(\theta_n)(I_n-\mathbb{E}I_n)^2\Big) 1_{A_n^c}\Big|\\
& \leq & \Big|F'(\mathbb{E}I_n)
\mathbb{E}(I_n-\mathbb{E}I_n)1_{A_n}\Big|+\frac{1}{2}\Big|\mathbb{E}F''(\theta_n)(I_n-\mathbb{E}I_n)^21_{A_n^c}
\Big|=\gamma_{1,n}+\gamma_{2,n},
\end{eqnarray*}
where $\theta_n \in [\mathbb{E}I_n-\kappa {\rm
Var}\,I_n,\mathbb{E}I_n+\kappa {\rm Var}\,I_n]$. Consequently, by
Cauchy-Shwartz and Chebyshev's inequalities we obtain
\begin{eqnarray*}
\gamma_{1,n}
 & = &    |F'(\mathbb{E}I_n) \mathbb{E}(I_n-\mathbb{E}I_n)1_{A_n}|
 \leq |F'(\mathbb{E}I_n)|\sqrt{\mathbb{E}(I_n-\mathbb{E}I_n)^2}\sqrt{\mmp(A_n)}\\
 & \leq & |F'(\mathbb{E}I_n)|\sqrt{{\rm Var}\,I_n}\sqrt{\frac{{\rm Var}\,I_n}{(\kappa {\rm Var}\,I_n)^2}}=\frac{1}{\kappa}|F'(\mathbb{E}I_n)|,
\end{eqnarray*}
which is $o(h(n))$ by condition (C7).

Finally, an appeal to condition (C4) allows us to conclude that
\begin{eqnarray*}
S_{4} & \leq & \frac{1}{2}\sup_{x\geq \mathbb{E}I_n/2}|F''(x)|
{\rm Var}\,I_n,
\end{eqnarray*}
which is $o(h(n))$ by condition (C5). The proof is complete.

\end{proof}
\end{thm}
Theorem \ref{theorem_general} and Theorem \ref{main_theorem_1}
justify the Algorithm given in the Introduction.

\section{Applications.}
\subsection{Exchangeable coalescents.}
\subsubsection{Number of collisions in beta$(a,1)$-coalescents.}

Let $X_n$ be the number of collisions in beta$(a,1)$-coalescent,
$a>0$, restricted to the set $\{1,\ldots,n\}$. Many results
concerning the asymptotics of $\me X_n^k$, $k\in\mn$, are known
\cite{Delmas, DIMR1,DIMR2, GneYak, IMM, MeirMoon, Pan1} but we
partially derive them again just in order to show how our method
works.

% Pitman \cite{Pit} and
%Sagitov \cite{Sag} independently introduced coalescent processes
%with multiple collisions. These Markovian processes with state
%space ${\cal E}$, the set of all equivalence relations
%(partitions) on $\mn$, are characterized by a finite measure
%$\Lambda$ on $[0,1]$ and are, hence, also called
%$\Lambda$-coalescents. We consider beta$(a,1)$-coalescents, which
%means that $$\Lambda(dx)=ax^{a-1}1_{(0,1)}(x)dx,$$ for some $a>0$.
It is known (see, for instance, \cite[Section 7]{IksMoe2}) that
the sequence $\{X_n,n\in\mn\}$ satisfy the distributional equality
$$
X_1=0,  \ \ X_n\od 1+X_{I_n}, \ \ n\geq 2,
$$
where $I_n$ is a random variable with distribution
$$\mmp\{I_n=n-k\}=\frac{\frac{(2-a)\Gamma(a+k-1)}{\Gamma(a)\Gamma(k+2)}}{1-\frac{\Gamma(a+n-1)}{\Gamma(a)\Gamma(n+1)}},
\ \ k=1,\ldots,n-1, \ \ n\geq 2,$$ if $a\neq 2$, and
$$\mmp\{I_n=n-k\}= \frac{1}{(h_n-1)(k+1)}, \ \ k=1,\ldots,n-1, \ \ n\geq
2,$$ where $h_n=\sum_{k=1}^n k^{-1}$, if $a=2$.

By Proposition \ref{is_unbounded} it follows that $\lin
\mathbb{E}X_n=+\infty$. It is also clear that no reduction of weights to
probabilities in the recurrence is needed.

\vspace{.2cm}\noindent {\sc Case} $0<a<1$ \cite{Delmas, GneYak,
IksMoe2}. Since
$$\mathbb{E}I_n=n-(1-a)^{-1}+o(1),$$ we can choose $$g(x)=x-{1\over 1-a} \ \ \text{and} \ \ h(x)= 1.$$
Then functional equation \eqref{iter_def} has an elementary
solution $g^\ast(x)=(1-a)x$. By Theorem \ref{theorem_general},
$\mathbb{E}X_n\sim g^{*}(n)\sim (1-a)n$.

\vspace{.2cm}\noindent {\sc Case} $a=1$ (Bolthausen-Sznitman
coalescent) \cite{DIMR1, IksMoe2, MeirMoon, Pan1}. Since
$$\mathbb{E}I_n=n-\log n+O(1),$$ we can choose $$g(x)=x-\log x \ \ \text{and} \ \ h(x)=
1.$$ From the relation $\frac{x}{\log x}=1+o(1)+\frac{x-\log
x}{\log (x-\log x)}$ and Theorem \ref{f_equiv} it follows that
${\rm Iter} (h,g,2)(x) \sim \frac{x}{\log x}$. An application of
Theorem \ref{main_theorem_1}\footnote{The only thing which may
require verification is condition $C3$. In the present situation,
$\mathbb{E}I_n-g(n)=O(1)$, and the derivative of $F(x)=x/\log x$
tends to zero, as $x\to\infty$. Therefore, condition $C3$ follows
by application of the mean value theorem.} gives
$\mathbb{E}X_n\sim \frac{n}{\log n}$.

\vspace{.2cm}\noindent {\sc Case} $a=2$ \cite{IMM}. Since
$$\mathbb{E} I_n=n-{n\over \log
n}+O\bigg(\frac{n}{\log^2 n}\bigg),$$ we can choose
$$g(x)=\bigg(x-{x\over \log x}\bigg)\ind_{(e,\infty)}(x) \ \ \text{and} \ \
h(x)=1.$$ From the relation $\frac{1}{2}\log^2
x=1+o(1)+\frac{1}{2}\log^2\Big(x-\frac{x}{\log x}\Big)$ and
Theorem \ref{f_equiv} we conclude that ${\rm Iter}(h,g,2)(x)\sim
\frac{1}{2}\log^2 x$. Direct calculations show that
$$\log^2 \mathbb{E}I_n -\log^2 g(n)=O\bigg({1\over \log n}\bigg)$$ and
$$\lin {1\over 2}\Big(\mathbb{E}\log^2 I_n-\log^2 \mathbb{E}I_n\Big)=
1 - \frac{\pi^2}{6},$$ which, in view of implication
\eqref{conclusion2}, yields $\mathbb{E}X_n\sim
\frac{3}{\pi^2}\log^2 n$.

\subsubsection{Functionals acting on the Poisson-Dirichlet coalescent.}\label{PDC}

Unlike the beta coalescents, the asymptotics of the moments of the
number of collisions $X_n$ in the Poisson-Dirichlet coalescent
does not seem to have been known so far (we refer to
\cite{Moehle2009, Sag2} for extensive information about this
particular coalescent with simultaneous multiple collisions).
Recall that this fact has served as an initial motivation for
developing the method reported in this article.

%\citeSchweinsberg \cite{Schweinsberg} characterizes exchangeable
%coalescents via a finite measure $\Xi$ on the infinite simplex
%$\Delta:=\{x=(x_1,x_2,\ldots):x_1\geq x_2\geq \cdots \geq
%0,\sum_{i=1}^{\infty}x_i\leq 1\}$. If $\Xi$ has no atom at zero
%and if $\Xi$ is concentrated on the subset $\Delta^*$ of points
%$x\in\Delta$ satisfying  $|x|=1$ then each
%$(k_1,\ldots,k_j)$-collision ($k_1,\ldots,k_j\in\mn$ with $k_1\geq
%\cdots k_j$ and $k_1\geq 2$) is occuring at the rate (see
%\cite[Eq. (11)]{Schweinsberg}
%\begin{eqnarray}\nonumber
%\psi_j(k_1,\ldots,k_j)=\int_{\Delta}\sumalldist x_1^{k_1}\cdots x_j^{k_j}\frac{\Xi(dx)}{(x,x)}.
%\end{eqnarray}
%An important and interesting example of exchangeable coalescent is
%the Poisson-Dirichlet coalescent with parameter $\theta>0$ (see ).
%In this case the characterizing measure $\Xi$ has (by definition)
%density $x\mapsto (x,x)$ with respect to the Poisson-Dirichlet
%distribution $\Pi_{\theta}$ with parameter $\theta>0$. Note that
%the measure $\Pi_{\theta}$ is concentrated on $\Delta^*$. From the
%calculations of Kingman (see \cite[Section 9.5]{KingmanPP} it
%follows that the Poisson-Dirichlet coalescents has rates
%\begin{eqnarray}\label{PD_freq}
%\psi_j(k_1,\ldots,k_j)=\frac{\theta^j}{[\theta]_k}\prod_{i=1}^{j}(k_i-1)!,
%\end{eqnarray}
%where $j,k_1,\ldots,k_j\in\mn$ with $k:=k_1+\ldots+k_j>j$, where $[\theta]_k:=\theta(\theta+1)\cdots (\theta+k-1)$.

It can be checked that the sequence $\{X_n, n\in\mn\}$ satisfies
the distributional equality
$$
X_1=0, \ \ X_n\od 1+X_{I_n}, \ \ n\geq 2,
$$
where
$$\mmp\{I_n=k\}=\frac{\theta^k}{[\theta]_n-\theta^n}s(n,k), \ \ k=1,\ldots,n-1, \ \ n\geq 2,$$ $s(n,k)$ is the unsigned
Stirling number of the first kind, and
$[\theta]_n=\theta(\theta+1)\ldots(\theta+n-1)$. This implies that
$$\mathbb{E}I_n=\theta \log n + O(1),\;\;{\rm Var}\,I_n=\theta
\log n +O(1).$$ Set $$g(x):=\theta\log x, \ \ h(x)=1 \ \
\text{and} \ \ g^{*}(x)={\rm Iter}(h,g,x_0)$$ for some
$x_0>\exp(2\theta\vee 1)$. Notice that the $g^{*}$ is a
generalized log-star function, which can be defined via functional
equation
$$
g^{*}(x)=1+g^{*}(\theta\log x),\;x>x_0.
$$
Let $F$ be a twice differentiable modification of the $g^{*}$ of
the form
$$
F(x)=
\left
\{
\begin{array}{ll}
1+F(\theta\log x),\;\;x>x_0,\\
\alpha_1 x^3+\alpha_2 x^2+\alpha_3 x, x\in [0,x_0],
\end{array}
\right.
$$
for some constants $\alpha_1,\alpha_2,\alpha_3$. It follows that,
for every fixed $j\in\mn$,
$$F'(x)=o\bigg(\frac{1}{x\log x \cdots \log^{\circ (j)}(x)}\bigg) \ \
\text{and} \ \ F''(x)=o\bigg(\frac{1}{x^2(\log x)^2 \cdots
(\log^{\circ (j)}(x))^2}\bigg).$$ An application of Theorem
\ref{main_theorem_1} yields $\mathbb{E}X_n\sim g^{*}(n)\sim F(n)$.
Analogously, we obtain
$$
\mathbb{E}X_n^k\sim (g^{*}(n))^k, \ \ k\in\mn.
$$

Other important functionals acting on the Poisson-Dirichlet
coalescent are the absorption time $T_n$ (in the biological
context, $T_n$ is the time back to the most recent common ancestor
of a sample of size $n$) and the total branch length $L_n$ of the
coalescent tree.

Corresponding distributional recurrences are
\begin{eqnarray*}
&   &   T_1=0,\;\;T_n\od  \tau_n+T_{I_n},\;\;n\geq 2,\\
&   &   L_1=0,\;\;L_n\od n\tau_n+L_{I_n},\;\;n\geq 2,
\end{eqnarray*}
where $\tau_n$ is a random variable with the exponential law with
parameter $g_n=1-\frac{\theta^n}{[\theta]_n}$ which is independent
of everything else.

Using induction on $k$, the fact that $\lin g_n= 1$ and Theorem
\ref{seq_equiv} we conclude that
$$\mathbb{E}T_n^k\sim (g^{*}(n))^k, \ \ k\in\mn.$$
An application of Chebyshev's inequality immediately leads to the
following weak laws of large numbers.
\begin{thm}
As $n\to\infty$,
\begin{eqnarray*}
\frac{X_n}{g^{*}(n)}\tp 1 \ \ \text{{\rm and }} \ \
\frac{T_n}{g^{*}(n)}\tp 1.
\end{eqnarray*}
\end{thm}

As far as $L_n$ is concerned, we can prove that
$$
\mathbb{E}L_n^k\sim k!n^k, \ \ k\in\mn.
$$
By the method of moments this immediately gives the following weak
convergence result.
\begin{assertion}
As $n\to\infty$,
\begin{eqnarray*}
\frac{L_n}{n}\dod L,
\end{eqnarray*}
where $L$ is a random variable with the standard exponential law.
\end{assertion}
In a recent preprint \cite{Moehle2009} the same result was
obtained by a different method. We thus omit further details.

\subsection{Examples from the analysis of algorithms.}
We will give new proofs of the results from \cite{NeiRus},\cite{Mahmoud} and \cite{RosQu}, respectively, by using our method.
\subsubsection{The Quickselect algorithm.}
Let $X_n$ be the number of comparisons that the Quickselect
algorithm needs to find $\min(x_1,\ldots,x_n)$ of a sample
$x_1,\ldots,x_n$. Then
$$
X_1=0,\;\;X_n\od n-1+X_{I_n}, \ \  n\geq 2,
$$
where $I_n=J_n\vee 1$, and $J_n$ is uniformly distributed on
$\{0,\ldots,n-1\}$. Since
$$\mathbb{E}I_n={n-1\over 2}+{1\over n},$$ we can choose $$g(x)={x+1\over 2} \ \
\text{and} \ \ h(x)=x-1.$$ Then functional equation
\eqref{iter_def} has elementary solutions $g^{*}(x)=2x+c$,
$c\in\mathbb{R}$. By Theorem \ref{theorem_general},
$\mathbb{E}X_n\sim g^{*}(n)\sim 2n$.
\subsubsection{The depth of a random node in a random binary search tree.}
The corresponding recurrence is
$$
X_0=-1, \ X_1=0, \ \ X_n\od 1+X_{I_n}, \ \ n\geq 2,
$$
where $\mmp\{I_n=k\}=2k/n^2$ for $k\in\{1,\ldots,n-1\}$ and
$\mmp\{I_n=0\}=1/n$. Since
$$\mathbb{E}I_n={(n-1)(2n-1)\over 3n},$$ we can choose $$g(x)=2x/3 \ \ \text{and} \ \ h(x)=1.$$
According to Example 1, the corresponding iterative function is
$$g^*(x)=\lceil \log_{{3\over 2}} x\rceil, \ \ x>1.$$
Since $\lin (\mathbb{E}\log^+ I_n-\log n)= -1/2$, it follows that
$\lin (\mathbb{E}f(I_n)-f(\mathbb{E}I_n))= 1-\frac{1}{2\log
(3/2)}$, where $$f(x)={\log^+ x\over \log (3/2)}.$$ Since
$f(x)\sim g^\ast(x)$ then, according to the Algorithm,
$\mathbb{E}X_n\sim 2\log (3/2)f(n)\sim 2\log n$.
\subsubsection{The Quicksort algorithm.}
Let $X_n$ denote the random number of comparisons needed to sort a
list of length $n$ by the Quicksort.  Then $X_0=X_1=0$, and
$$
X_n\od n-1+ X_{I_n-1}+X'_{n-I_n}, \ \ n\geq 2,
$$
where $\{X'_n,n\in\mno\}$ is an independent copy of
$\{X_n,n\in\mno\}$ which is independent of $I_n$ having the
uniform distribution on $\{1,\ldots,n\}$. Set
$a_n:=\mathbb{E}X_n$, then $a_0=a_1=0$ and
$$
a_n=n-1+\sum_{k=0}^{n-1}\frac{2}{n}a_k, \ \ n\geq 2.
$$
The reduction of weights to probabilities can be made by the substitution
$a'_n:=a_n/(n+1)$ which yields
$$
a'_n=\frac{n-1}{n+1}+\sum_{k=0}^{n-1}\frac{2(k+1)}{n(n+1)}a'_k, \
\ n\geq 2.
$$
Using the same arguments as in the previous example we obtain $a'_n\sim 2\log n$. Therefore, 
$\mathbb{E}X_n\sim 2n\log n$, which is well known asymptotic for Quicksort.
\footnote{
First result concerning the complexity of the (non-randomized) Quicksort algorithm with $O(n\log n)$ asymptotic goes back to the pioneer work by Hoar \cite{Hoar}. For complete analysis of the Quicksort and its different modifications we refer to survey \cite{QuicksortSedg}.}

\subsection{Limitations of the method.}
\begin{enumerate}
\newcommand{\thenumi}{\alph{enumi}}
\renewcommand{\labelenumi}{(\thenumi)}
\item An indispensable requirement of our method to work is the divergence of $a_n$, the solution to (\ref{rec_mean}).
In particular, our method cannot detect the convergence of $a_n$
to a constant.
\item It may be difficult to guess which elementary function has the same asymptotics as a given iterative function.
\item If condition (\ref{conclusion2}) holds for some $c\neq 0$, it may be hard to calculate the constant $c$ explicitly.
Therefore, it seems that a natural assumption for the method to
work is (\ref{conclusion1}) rather than (\ref{conclusion2}).
Condition (\ref{conclusion1}) holds if the solution is nearly
linear and the variance of index $I_n$ grows not too fast (precise
statements are made in Theorem \ref{main_theorem_1}). For
instance, the mean number of collisions in the Bolthausen-Sznitman
and Poisson-Dirichlet coalescents exhibit the asymptotic behavior
of this type.
\end{enumerate}

\section{Appendix.}
\subsection{Some properties of iterative functions.}\label{mod}
For the given strictly increasing continuous
function $g$, there exists the unique inverse function $g^{-1}$
which defines the sequence $\{A_n, n\in\mno\}$ as follows
\begin{eqnarray}\label{Adef}
A_0=0,\;\;A_i:=(g^{-1})^{\circ (i-1)}(x_0),\;\;i\in\mn.
\end{eqnarray}

\begin{lemma}\label{lem1}
Assume that $g,h,k\in C^{(m)}[x_0,+\infty)$ and $F={\rm
Iter}(h,g,x_0,k)$ is $m$-times differentiable at $x_0$. Then $F$
is $m$-times differentiable on $[x_0,+\infty)$.
\begin{proof}
We only treat the case $m=1$, as, for $m=2,3,\ldots$, the proof is
the same. Since $F$ is a sum of compositions of
$C^{(1)}[x_0,+\infty)$ functions, it is differentiable on
$[x_0,+\infty)\backslash \{A_i, i\in\mn\}$. Therefore, we only have to check the continuity and differentiability at points $\{A_i,\;i\in\mn\}$.\\
\textit{First step.} Proof of continuity. By the assumption, $F$
is continuous at $A_1=x_0$, i.e.,
\begin{eqnarray}\label{continuity}
k(x_0)=h(x_0)+k(g(x_0)).
\end{eqnarray}
For fixed $k\geq 2$, we have from (\ref{sum_representation})
\begin{eqnarray}\label{lem1f1}
F(A_k-0)
 & = &\sum_{i=1}^{k-1}h(g^{\circ (i-1)}(A_k-0))+k(g^{\circ(k-1)}(A_k-0)),
\end{eqnarray}
and
\begin{eqnarray}\label{lem1f2}
F(A_k+0)
 & = &\sum_{i=1}^{k}h(g^{\circ (i-1)}(A_k+0))+k(g^{\circ(k)}(A_k+0)).
\end{eqnarray}
Use now (\ref{continuity}) and continuity of $h$ and $g$ to obtain
\begin{eqnarray*}
&   &\hspace{-2.4cm}
   F(A_k+0)-F(A_k-0)\\
   \hspace{-1cm}& = & h(g^{\circ (k-1)}(A_k))+k(g^{\circ(k)}(A_k))-k(g^{\circ (k-1)}(A_k))\\
   & = & h(x_0)+k(g(x_0))-k(x_0)\overset{(\ref{continuity})}{=}0.
\end{eqnarray*}
\textit{Second step.} Proof of differentiability. The
differentiability of $F$ at $x_0$ implies that
\begin{eqnarray}\label{differentiable}
k'(x_0)=h'(x_0)+k'(g(x_0))g'(x_0).
\end{eqnarray}
For $k\geq 2$, using (\ref{lem1f1}) and (\ref{lem1f2}) yields
\begin{eqnarray*}
 & & F_{-}'(A_k)=\lim_{x\to A_k-0}{d\over dx}\bigg(\sum_{i=1}^{k-1}h(g^{\circ (i-1)}(x))+k(g^{\circ
(k-1)}(x))\bigg),\\
 & & F_{+}'(A_k)=\lim_{x\to A_k+0}{d \over dx}\bigg(\sum_{i=1}^{k}h(g^{\circ (i-1)}(x))+k(g^{\circ(k)}(x))\bigg).
\end{eqnarray*}
Consequently,
\begin{eqnarray*}
   &   & \hspace{-1cm}
F_{+}'(A_k)-F_{-}'(A_k)\\
 & = &\lim_{x\to A_k+0}\frac{d}{dx}h(g^{\circ(k-1)}(x))+k(g^{\circ (k)}(x))-\lim_{x\to A_k-0}\frac{d}{dx}k(g^{\circ (k-1)}(x)).
\end{eqnarray*}
Set $u(x):=g^{\circ (k-1)}(x)$, then $u(A_k+0)=u(A_k-0)=u(A_k)=x_0$ and
\begin{eqnarray*}
   &   & \hspace{-1cm}
         F_{+}'(A_k)-F_{-}'(A_k)\\
   & = & \lim_{x\to A_k+0}\frac{d}{dx}h(u(x))+k(g(u(x)))-\lim_{x\to A_k-0}\frac{d}{dx}k(u(x))\\
   & = & \lim_{x\to A_k+0}(h'(u(x))+k'(g(u(x)))g'(u(x)))u'(x)-\lim_{x\to A_k-0}k'(u(x)))u'(x)\\
   & = & (h'(x_0)+k'(g(x_0))g'(x_0)-k'(x_0))u'(x_0)=0,
\end{eqnarray*}
by \eqref{differentiable}. The proof is complete.
\end{proof}
\end{lemma}
From this lemma it follows that the function $F$ is $m$-times
differentiable provided it satisfies conditions
\begin{eqnarray}\label{k_conds}
\begin{array}{ll}
k(x_0)=h(x_0)+k(g(x_0)),\\
k'(x_0)=h'(x_0)+k'(g(x_0))g'(x_0),\\
\ldots\ldots\ldots\ldots\ldots\ldots\ldots\ldots\ldots\\
k^{(m)}(x_0)=h^{(m)}(x_0)+(k(g(x_0)))^{(m)}.
\end{array}
\end{eqnarray}
The following lemma proves the existence of such a function
$k(x)$.
\begin{lemma}\label{lem2}
Assume that
$W\Big(x-g(x),\ldots,x^{m+1}-g^{m+1}(x)\Big)\Big|_{x=x_0}\neq 0$.
Then there exists a function $k(x)=\sum_{i=1}^{m+1}\alpha_ix^i$
which satisfies (\ref{k_conds}).
\begin{proof}
Plugging the representation $k(x)=\sum_{i=1}^{m+1}\alpha_ix^i$
into (\ref{k_conds}) gives the system of linear equations
$$
\begin{array}{ll}
\Big(\alpha_1(x_0-g(x_0))+\ldots+\alpha_{m+1}(x_0^{m+1}-g^{m+1}(x_0))\Big)=h(x_0),\\
\Big(\alpha_1\frac{d}{dx}(x-g(x))+\ldots+\alpha_{m+1}\frac{d}{dx}(x^{1+m}-g^{m+1}(x))\Big)\Big|_{x=x_0}=h'(x_0),\\
\ldots\ldots\ldots\ldots\ldots\ldots\ldots\ldots\ldots\\
\Big(\alpha_1\frac{d^m}{dx^m}(x-g(x))+\ldots+\alpha_{m+1}\frac{d^m}{dx^m}(x^{m+1}-g^{m+1}(x))\Big)\Big|_{x=x_0}=h^{(m)}(x_0).
\end{array}
$$
The determinant of this system is  $W(x_0)$ which is not equal to
zero by the assumption. Therefore, the system has a unique
solution which implies that the function $k$ is well defined and
satisfies conditions (\ref{k_conds}).
\end{proof}
\end{lemma}
\subsection{Inhomogeneous terms of recursion (\ref{rec_mean}) and iterative functions.}
\begin{thm}\label{seq_equiv}
Suppose that $\{a_n,n\in\mn\}$ and $\{a'_n,n\in\mn\}$ satisfy the
recurrences
\begin{eqnarray}\label{rec1}
a_n=b_n+\sum_{k=1}^{n-1}p_{nk}a_k,n\geq N,
\end{eqnarray}
and
\begin{eqnarray}\label{rec2}
a'_n=b'_n+\sum_{k=1}^{n-1}p_{nk}a'_k,n\geq N,
\end{eqnarray}
respectively. Suppose that $b_n\geq 0$ for $n\geq N$ and $\lin
a_n=+\infty$. Then
\newcounter{count1}
\setcounter{count1}{1}
\begin{list}{\Roman{count1}.}%
{\usecounter{count1}}
\item  $b'_n\sim b_n$, $n\to\infty$ implies $a'_n\sim a_n$, $n\to\infty$, and
\item  $b'_n=o(b_n)$, $n\to\infty$ implies $a'_n=o(a_n)$, $n\to\infty$.
\end{list}
\begin{proof}[Proof of (I)]
We exploit the idea of proof of \cite[Proposition
3]{gnedin_sieve}. Suppose there exists $\varepsilon_0>0$ such that
$a_n>(1+\varepsilon_0)a'_n$ for infinitely many $n$. Since $\lin
a_n=+\infty$, we can pick $\varepsilon\in (0,\varepsilon_0]$ such
that for any $c>0$ the inequality $a_n>(1+\varepsilon)a'_n+c$
holds for infinitely many $n$. Let $n_c$ be the minimal such $n$.
Since $\underset{c\to\infty}{\lim}n_c=+\infty$, without loss of
generality we can assume that $n_c>N$. For $n\leq n_c-1$ we have
$a_n<(1+\varepsilon)a'_n+c$ which implies
$$
(1+\varepsilon)a'_{n_c}+c<a_{n_c}=b_{n_c}+\sum_{k=1}^{n_c-1}p_{n_c
k}a_k<b_{n_c}+c+(1+\varepsilon)\sum_{k=1}^{n_c-1}p_{n_ck}a'_k.
$$
Simplifying the last expression gives
$1+\varepsilon<b_{n_c}/b'_{n_c}$. Sending $c\to\infty$ leads to
$\varepsilon<0$, which is a contradiction. Thus we have proved
that
$$
\limsup_{n\to\infty} \frac{a_n}{a'_n}\leq 1.
$$
A symmetric argument proves the converse inequality for $\liminf$.

\textit{Proof of (II)} proceeds by applying the already
established part (I) to the sequences $\{a_n,n\in\mn\}$ and
$\{a_n-a'_n,n\in\mn\}$ and noting that the relation $b_n\sim
b_n-b'_n$ implies $a_n\sim a_n-a'_n$. The proof is complete.
\end{proof}
\end{thm}
Using a similar reasoning one can prove the following.
\begin{thm}\label{f_equiv}
Let the triples $(h_1,g,x_0)$ and $(h_2,g,x_0)$ generate the
iterative functions $f_1$ and $f_2$, respectively. %So they are
Assume that $\lix f_1(x)=+\infty$. Then
\newcounter{count2}
\setcounter{count2}{1}
\begin{list}{\Roman{count2}.}%
{\usecounter{count2}}
\item  $h_2(x)\sim h_1(x)$, $x\to\infty$ implies $f_2(x)\sim f_1(x)$, $x\to\infty$, and
\item  $h_2(x)=o(h_1(x))$, $x\to\infty$ implies $f_2(x)=o(f_1(x))$, $x\to\infty$.
\end{list}

\end{thm}
\subsection{Sufficient condition for the divergence of solutions to (\ref{rec_mean}).}
A simple sufficient condition for $\lin a_n= +\infty$ is given
next.
\begin{assertion}\label{is_unbounded}
Assume that the sequence $\{a_n,n\in\mn\}$ satisfies
(\ref{rec_mean}). If $I_n\tp \infty$ and
$\underset{n\to\infty}{\lim \inf}\, b_n=b>0$, then $\lin
a_n=+\infty$.
\begin{proof}
From recurrence (\ref{rec_mean}) we obtain
\begin{eqnarray*}
 a_n
 & = & b_n+\sum_{k=1}^{n-1}p_{nk}a_k=b_n+\sum_{k=1}^{M-1}p_{nk}a_k+\sum_{k=M}^{n-1}p_{nk}a_k\\
 & \geq & b_n+\Big(\inf_{1\leq k< M}a_k\Big)\sum_{k=1}^{M-1}p_{nk}+\Big(\inf_{M\leq k\leq n-1}a_k\Big)\sum_{k=M}^{n-1}p_{nk}.
\end{eqnarray*}
Sending $n\to\infty$ gives $\underset{n\to\infty}{\lim
\inf}\,a_n\geq b+\underset{k\geq M}{\inf}\, a_k$. Letting
$M\to\infty$ leads to $\underset{n\to\infty}{\lim \inf}\, a_n\geq
b+\underset{n\to\infty}{\lim \inf}\,a_n$ which completes the
proof.
\end{proof}
\end{assertion}

\vskip0.1cm \noindent {\bf Acknowledgement}. I would like to thank
my scientific advisor Alexander Iksanov for introducing me into
the subject, careful reading of the manuscript and numerous
helpful comments. Thanks are also due to Martin M\"{o}hle for
drawing our attention to the problem about asymptotics of the
Poisson-Dirichlet coalescent. Finally I am grateful to the referee 
for careful reading of the manuscript.


\begin{thebibliography}{99}
\bibitem{Bruhn}{\sc Bruhn, V.} (1996) Eine Methode zur asymptotischen Behandlungeiner Klasse von Rekursionsgleichungen mit einer Anwendung in der stochastischen Analyse des Quicksort-Algorithmus. {\em Dissertation, Christian-Albrechts-Universit\"at zu Kiel.}

\bibitem{Cutsem}{\sc van Cutsem, B. and Ycart, B.} (1994). Renewal-type
behaviour of absorption times in Markov chains. {\em Adv. Appl.
Prob.} {\bf 26,} 988--1005.

\bibitem{Delmas}{\sc Delmas, J.-F., Dhersin, J.-S. and Siri-Jegousse, A.}
(2008). Asymptotic results on the length of coalescent trees. {\em
Ann. Appl. Probab.} {\bf 18,} 997Ц-1025.

\bibitem{Devillers}{\sc Devillers, O.} (1992). Randomization yields simple $O(n\log^* n)$ algorithms for difficult $\Omega(n)$ problems. {\em Internat. J. Comput. Geom. Appl.} {\bf 2}, 621--635.

\bibitem{Drmota}{\sc Drmota, M.} (2009). {\em Random trees: An interplay between combinatorics and probability}, Springer.


\bibitem{DIMR1}{\sc Drmota, M., Iksanov, A., Moehle, M. and
Roesler, U.} (2007). Asym\-ptotic results concerning the total
branch length of the Bolthausen-Sznitman coalescent.
{\em Stoch. Process. Appl.} {\bf 117}, 1404--1421.

\bibitem{DIMR2}{\sc Drmota, M., Iksanov, A., Moehle, M. and
Roesler, U.} (2009). A limiting distribution for the number of
cuts needed to isolate the root of a random recursive tree.
{\em Random Struct. Algorithms} {\bf 34}, 319--336.

\bibitem{FlajSedg}{\sc Flajolet, P. and Sedgewick, R.} (2008). {\em Analytic combinatorics}, Cambridge University Press.

\bibitem{gnedin_sieve}{\sc Gnedin, A.~V.} (2004). The Bernoulli sieve.
{\em Bernoulli} {\bf 10,} 79--96.

\bibitem{GneIksMoe}{\sc Gnedin, A., Iksanov, A. and M\"ohle, M.} (2008).
On asymptotics of exchangeable coalescents with multiple
collisions. {\em J. Appl. Probab.} {\bf 45,} 1186--1195.

\bibitem{GINR}{\sc Gnedin, A., Iksanov, A., Negadajlov, P. and Roesler, U.} (2009).
The Bernoulli sieve revisited. {\em Ann.\,Appl.\,Prob.} {\bf 19},
1634--1655.

\bibitem{GneYak}{\sc Gnedin, A. and Yakubovich, Y.} (2007). On the
number of collisions in $\Lambda$-coalescents.
{\em Electron. J. Probab.} {\bf 12}, 1547--1567.

\bibitem{Haas} {\sc Haas, B. and Miermont, G.} (2009). Self-similar scaling limits of non-increasing Markov chains, preprint available
at www.arXiv.org.

\bibitem{Hoar} {\sc Hoar, C.~R.} (1962). Quicksort. {\em The Computer Journal} {\bf 5(1),} 10-16.

\bibitem{KnutGrin}{\sc Greene D.H. and Knuth D.E.} (1990). {\em Mathematics for the analysis of algorithms} 3d Edition, Burkhauser.

\bibitem{IMM}{\sc Iksanov, A., Marynych, A. and M\"ohle, M.} (2009).
On the number of collisions in beta(2,b)-coalescents. {\em
Bernoulli} {\bf 15,} 829--845.

\bibitem{IksMoe2}{\sc Iksanov, A. and M\"ohle, M.} (2008). On the
number of jumps of random walks with a barrier. {\em Adv. Appl.
Probab.} {\bf 40}, 206--228.

\bibitem{IksNegad}{\sc Iksanov, A. and Negadajlov, P.} (2008). On
the number of zero increments of a random walk with a barrier.
{\em Discrete Math. Theor. Comput. Sci.} Proceedings Series Volume
AI, 247--254.

\bibitem{Karp}{\sc Karp R.M.} (1994). Probabilistic recurrence relations.
{\em Journal of the Association for Computer Machinery} {\bf 41,}
1136--1150.

\bibitem{Kesten}{\sc Kesten, H.} (1980). The number of distinguishable alleles according to the Ohta-Kimura model
of neutral mutation. {\em J. Math. Biol.} {\bf 10,} 167--187.

%\bibitem{KingmanPP}{\sc Kingman, J.F.C.} (1993). {\em Poisson processes,} Oxford University Press.

\bibitem{Mahmoud}{\sc Mahmoud H.} (1992). {\em Evolution of random search trees}, Wiley, New York.

\bibitem{MeirMoon}{\sc Meir, A. and Moon J.W.} (1974). Cutting down recursive trees.
Ц {\em Math. Biosci.} {\bf 21,}  173-Ц181.

\bibitem{Moehle2009}{\sc M\"ohle, M.} (2009). Coalescent processes
without proper frequencies and applications to the two-parameter
Poisson-Dirichlet coalescent, preprint available at
http://www.mathematik.uni-tuebingen.de/$\sim$ moehle/moehle.html


\bibitem{NeiRus}{\sc Neininger, R. and R\"uschendorf, L.} (2004).
On the contraction method with degenerate limit equation. {\em
Ann. Probab.} {\bf 32,} 2838--2856.

\bibitem{Pan1}{\sc Panholzer, A.} (2004). Destruction of recursive trees. In: {\em
Mathematics and Computer Science III}, Birkh\"{a}user, Basel,
267--280.

\bibitem{Pan2}{\sc Panholzer, A.} (2006). Cutting down very simple trees. {\em
Quaest. Math.} {\bf 29,} 211--227.

%\bibitem{Pit}{\sc Pitman, J.} (1999). Coalescents with multiple
%collisions. {\em Ann. Probab.} {\bf 27,} 1870--1902.

\bibitem{RosQu}{\sc R\"{o}sler, U.} (1991). A limit theorem for "Quicksort". {\em
RAIRO, Inform. Theor. Appl.} {\bf 25,} 85--100.

\bibitem{Roesler99}{\sc R\"{o}sler, U.} (2001). On the analysis of stochastic divide and conquer algorithms. {\em Algorithmica} {\bf 29,} 238--261.

%\bibitem{Sag}{\sc Sagitov, S.} (1999). The general coalescent with asynchronous
%mergers of ancestral lines. {\em J. Appl. Prob.} {\bf 36,} 1116--1125.

\bibitem{Sag2}{\sc Sagitov, S.} (2003). Convergence to the coalescent with simultaneous multiple merges. {\em J. Appl. Prob.} {\bf 40,} 839--854.

\bibitem{QuicksortSedg}{\sc Sedgewick, R.} (1977). The analysis of Quicksort Programs. {\em Acta Informatika} {\bf 7(4),} 327--355.

%\bibitem{Schweinsberg}{\sc Schweinsberg J.} (2000). Coalescents with simultaneous multiple  collisions. {\em Electron. J. Probab.} {\bf 5,} 1--50.

\end{thebibliography}
\end{document}